\documentclass[runningheads]{llncs}
\usepackage[T1]{fontenc}
\usepackage{amsfonts,bm}
\usepackage{amssymb,amsmath}
\usepackage{graphicx,subfigure}
\usepackage{cite}
\usepackage{comment}
\usepackage{algorithm2e}
\usepackage{algpseudocode}
\RestyleAlgo{ruled}
\usepackage{ifthen}     

\renewcommand{\epsilon}{\varepsilon}
\renewcommand{\phi}{\varphi}
\renewcommand{\d}{\,\mathrm{d}}

\renewcommand{\epsilon}{\varepsilon}
\renewcommand{\phi}{\varphi}
\renewcommand{\d}{\,d}

\newcommand{\R}{\mathbb{R}}

\begin{document}
\title{On Minimum-Dispersion Control of Nonlinear Diffusion Processes 
}
%
%
\author{Roman Chertovskih\inst{1}\orcidID{0000-0002-5179-4344} \and
 Nikolay Pogodaev\inst{2}\orcidID{0000-0002-7062-1764} \and
Maxim Staritsyn\inst{1}\orcidID{0000-0003-3938-3128} \and   
A.~Pedro Aguiar\inst{1}\orcidID{0000-0001-7105-0505}}%
\authorrunning{R. Chertovskih {\it et al.}}
%
\institute{Research Center for Systems and Technologies (SYSTEC), ARISE, Department of Electrical
and Computer Engineering,
Faculdade de Engenharia, Universidade do Porto, Rua Dr. Roberto Frias, s/n 4200-465, Porto, Portugal\\
\email{\{roman,staritsyn,pedro.aguiar\}@fe.up.pt}\\
\and
Dipartimento di Matematica ``Tullio Levi-Civita'' (DM), University of Padova, Via Trieste, 63 - 35121 Padova, Italy\\
\email{nickpogo@gmail.com}}
\maketitle              
\begin{abstract}
This work collects some methodological insights for numerical solution of a ``minimum-dispersion'' control problem for nonlinear stochastic differential equations, a particular relaxation of the covariance steering task. 
The main ingredient of our approach is the theoretical foundation called $\infty$-order variational analysis. This framework consists in establishing an exact representation of the increment ($\infty$-order variation) of the objective functional using the duality, implied by the transformation of the nonlinear stochastic control problem to a linear deterministic control of the Fokker-Planck equation.
The resulting formula for the cost increment analytically represents a ``law-feedback'' control for the diffusion process. This control mechanism enables us to learn time-dependent coefficients for a predefined Markovian control structure using Monte Carlo simulations with a modest population of samples.
Numerical experiments prove the vitality of our approach.


\end{abstract}

\section{Introduction: Optimal Covariance Steering vs Minimum-Dispersion Control}

This study was motivated by our interest in the \emph{covariance steering problem} (CSP) \cite{Hotz,GRIGORIADIS1997569, ChenY,BAKOLAS201861,Okamoto,TSOLOVIKOS202159,yu2023covariance}, a prominent aspect of applied stochastic control theory. Essentially, CSP involves guiding the random state of a stochastic process to one with a predefined mean and covariance matrix at a fixed terminal time, given an initial Gaussian probability distribution. As an intuitive example, one may have in mind the challenge of safely landing an aircraft in a noisy environment within a designated ``safe zone'' with a reasonable probability. Notably, CSP can be viewed as a sort of stochastic optimal control under mass transportation constraints.

CSP is known to have a closed-form solution in the case of Gaussian initial distribution, linear dynamics, and linear-quadratic cost function. The result is derived by the standard algebraic arguments based on the solution of a certain Riccati equation, akin to those found in the framework of linear-quadratic regulators. Recent attempts to address the nonlinear case \cite{TSOLOVIKOS202159,yu2023covariance} basically use the linearization of the state dynamics and employ reasoning from the linear case.  

Naturally,  the presence of nonlinear dynamics requires a relaxation of the mean-covariance fitting constraint, assuming that uncertainty in the terminal state decreases or does not propagate.  In \cite{BAKOLAS201861,TSOLOVIKOS202159}, such a relaxation is reached by imposing a condition on the ``growth'' of the covariance matrix. Note, however, that even when starting with a Gaussian random variable, completely characterized by its mean and covariance, nonlinearity in the driving vector field typically disrupts the Gaussian structure. Consequently, second-order statistics become insufficient for capturing the shape of the resulting probability distribution, and the associated steering problem should incorporate higher-order moments.

In this paper, we explore a different yet closely related class of problems to be referred to as \emph{minimum-dispersion control}, which captures the idea of steering the mean of a stochastic population toward a predefined target while considering a suitable higher-order statistical measure of scatter around the mean. This can be considered as a specific form of penalizing the mass transportation constraint in CSP such that, even in the linear-Gaussian case, a satisfactory solution to our problem implies a solution to the corresponding relaxed CSP. 

Essentially, the problem involves nonlinear stochastic control with Markovian strategies and a specific cost function, which we tackle numerically. It is widely understood that, under generally accepted regularity assumptions, the \emph{nonlinear stochastic} problem can be reformulated as an optimization involving a linear functional of the law of the random state subject to the Fokker-Planck equation in the space of probability measures. Consequently, this problem becomes \emph{linear} in the state variable and  \emph{deterministic}, allowing for standard duality arguments. Leveraging this duality, we derive an exact representation of the increment ($\infty$-order variation) of the cost functional, from which we can extract the structure of the descent control in the form of feedback dependence on the law --- probability measure. This designed law-feedback control is then employed to synthesize a descent Markovian control of a predefined structure using the Krasovskii-Subbotin algorithm \cite{krasovskii2012control}. To implement this approach, we combine Monte Carlo approximation of the adjoint state --- a solution to the backward Kolmogorov equation (the dual of the Fokker-Planck PDE) --- with an empirical approximation of the law of the target random variable. We demonstrate the effectiveness and practical applicability of this approach through a series of numerical experiments.

\section{Optimal Control Problem}

On a fixed time interval $I \doteq [0,T]$, $T>0$, consider a standard optimal stochastic control problem in the Mayer's form:
\[
(P) \qquad   \min \mathcal I[u] =  \mathbb E \, \ell(X_T[u]), \quad u \in \mathcal U.
\]
Here, $\ell \colon \R^n \to \R$ is a given cost function, $u$ is a control, and $X[u]$ is the corresponding trajectory being a strong solution to the nonlinear stochastic differential equation (SDE)
\begin{equation}
   X_t = x_0 + \int_0^t f_\tau(X_s, u_s) \d s + \int_0^t \sigma_s(X_s, u_s) \d W_s, \quad t \in I.\label{SDE}    
\end{equation}
Precisely, $X[u]\colon (t, \omega) \mapsto X_t[u](\omega)$ is a random process $I \times \Omega \to \R^n$ on a complete filtered probability  space $(\Omega, \mathcal F, \mathcal F_t^{W, X_0}, \mathbb P)$ with the filtration $t \mapsto \mathcal F_t^{W, X_0}$ of the sigma-algebra $\mathcal F$ generated by an $m$-dimensional Wiener process $W$ with independent components and a given initial random state $X_0 \in (\Omega, \mathcal F, \mathbb P)$, independent of $W$. The second integral in \eqref{SDE} is understood in \^{I}to's sense. 

By $\mathbb E$, $\mathbb V$ and $\mathbb K$ we denote the expectation, variance and covariance matrix of a random vector $\Omega \to \R^n$, respectively. 

The class $\mathcal U$ of feasible control actions, and the regularity of the drift and diffusion coefficients $f$ and  $\sigma$ are discussed in the next paragraphs. 

\subsection{Control Inputs}

An imperative of the stochastic control theory is the use of so-called Markovian control strategies $w\colon I \times \R^{n} \to U$ which are, in general, measurable function of both time $t$ and spatial position $x \in \R^{n}$ to a given set $U \subset \R^{d}$. However, this choice poses two serious challenges: i) a lack of regularity in the resulting vector field, leading to the existence of only weak solutions to the SDE \eqref{SDE}, and ii) impractical complexity of its implementation in real-life scenarios.

In this work, we sidestep these issues by adhering to control strategies $w(t, x) = \sum_j \xi_j(x) u_j(t)$ of a predefined structure w.r.t. $x$, involving a finite collection of given sufficiently regular functions $\xi_j \colon \R^n \to \R^d$.\footnote{In practice, the family $\{\xi_j\}$ would be smooth and bracket-generating (i.e. satisfying the H\"{o}rmander's condition).} We thus reframe problem $(RP)$ in the class \(
\mathcal U \doteq L_\infty(I; U)
\) of \emph{ensemble control strategies}, i.e., measurable functions $I \to U$ of the \emph{time variable only}.

\subsection{Standing Assumptions}

We make the following standard \underline{regularity hypotheses}: i) $U$ is compact; ii) the functions 
${f} \colon I \times \R^n  \times U \to \R^n$, $f=f_t(x, \upsilon)$, and $\sigma\colon I  \times \R^n\times U\to \R^{m\times n}$, $\sigma = \sigma_t(x)$, are bounded, measurable in $t \in I$, continuous in $(x, \upsilon)\in \R^n \times U$, and satisfy the Lipschitz and sublinear growth conditions w.r.t. $x \in \R^n$ uniformly in $(t, \upsilon)$; iii) $\mathbb E\big[|X_0|^2\big] < \infty$, and iv) $\ell\in C^2(\R^n)$ and satisfies the quadratic growth condition.

Recall that assumptions ii) and iii) guarantee the existence of a square integrable strongly unique strong solution to the SDE \eqref{SDE} \cite[Thm.~1.3.15]{Pham2021}. 

\subsection{Fokker-Planck Control Framework}\label{sec:fp}

Recall that the law of a random vector $X \colon \Omega \to \R^n$, denoted by $\mu = {\rm Law}(X)$ or $X \sim \mu$, is a probability measure $\mu = \mu_X \in \mathcal P(\R^n)$ defined via its action on test-functions $\phi\in \bm C_0(\R^n)$ (continuous real functions vanishing at infinity) as 
\(
\displaystyle \int_{\R^n} \phi \d \mu = \mathbb E(\phi\circ X).
\)
As a preliminary step in our analysis, we convert the nonlinear stochastic control problem $(P)$ into an equivalent state-linear deterministic optimization problem concerning the law dynamics $t \mapsto \mu_t = \text{Law}(X_t)$:
\begin{align}
(RP) \quad \min \quad & \mathcal J[u]= \int_{\R^d} \ell \, \d \mu_T[u], \ \ u \in \mathcal U;\label{I-func}\\
\mbox{subject to} \ \ &\partial_t\mu = \mathcal L_t^*(u_t) \, \mu, \ \,  t \in I,\mbox{ and } \label{PDE}  
X_0 \sim \mu_0,
\end{align}
where $\mathcal L^*_t(\upsilon)$ stands for the formal adjoint of the second-order elliptic operator $\mathcal L_t(\upsilon)$, $(t,\upsilon) \in I \times U$, acting as
\[
\mathcal L_t(\upsilon) \, \phi \doteq  f_t(x, \upsilon) \cdot \nabla_{x} \phi + {\rm Tr}\left(\nabla^2_{xx} \phi \,  D_t(x, \upsilon)\right) \quad \forall \phi \in C^2(\R^n).
\]
Hereinafter, ${\rm Tr}$ denotes the trace of a matrix, and $D \doteq \frac{1}{2} \sigma^{\rm T} \sigma$. The operator family $(t, \upsilon) \mapsto \mathcal{L}_t(\upsilon)$ is referred to as the generating family for controlled \^{I}to's diffusion \eqref{SDE}. 

The PDE \eqref{PDE} is understood in the following sense: for any $\phi \in C^\infty_c(\R^n)$ (smooth functions with compact support), and all $\tau, t \in I$, $\tau \leq t$, it holds
\[
 \int_{\R^n}\phi \d \mu_t - \int_{\R^n} \phi \d\mu_\tau = \int_\tau^t  \d s \int_{\R^n}\mathcal L_s(u_s) \phi  \d \mu_s,  
\]
and the initial condition means that $\lim\limits_{t \to 0}\mu_t = \mu_0$ in the corresponding weak* topology. We refer to \cite[\S~9]{bogachev2015fokker} for conditions ensuring the uniqueness of a (measure-valued) solution to \eqref{PDE}. 

It is noteworty that the transformation from the original problem $(P)$ to its relaxed counterpart $(RP)$ has become a modern imperative in stochastic optimal control theory. This transformation is widely utilized in computational contexts, as evidenced by numerous studies such as \cite{Annunziato2013, Roy2016, Roy2018, Breitenbach2020, Fleig2016, Fleig2017, Borzi2011}.

Some results on the existence of a solution to the optimization problems $(P)$ and $(RP)$ can be found in \cite[Thm. 6.3]{FlemingRishel} and  \cite[Thm.~4]{Anita2021}.

\subsection{Some Cost Functionals Measuring Dispersion of Samples}\label{sec:co}

Let us delve into some specifications of problem $(P)$ within the \emph{minimum-dispersion control} framework outlined in the Introduction. 

One natural approach to capture the idea of ``minimizing the scatter of random terminal states around the target mean $\hat x \in \mathbb{R}^n$ with regards to \emph{higher-order statistics}'' is to define the cost function $\ell$ as a sum of mixed central moments:
\[
\mathfrak{m}_{\alpha}(\mu; \hat x) \equiv \int_{\mathbb{R}^n} \prod_{j=1}^{\rm p}(x_j - \hat x_j)^{\alpha_j} \, \mathrm{d} \mu(x),
\]
where $\alpha=(\alpha_1, \alpha_2, \ldots, \alpha_{\rm p})$, $\alpha_j \in \mathbb{N}$, $j = 1, \ldots, {\rm p}$, is a multi-index of order $|\alpha| \equiv \sum_{j=1}^{\rm p} \alpha_{j}$, and $x_j$ represents the $j$th component of $x \in \mathbb{R}^n$. The corresponding functional \eqref{I-func} penalizes the $|\alpha|$-deviation of the random vector $X_T[u]$ from the target $\hat x$. Notably, the second-order characteristic:
\(
\ell(x) = 
\int_{\mathbb{R}^n} \|x - \hat x\|^2 \, \mathrm{d} \mu(x),
\)
can still serve as a suitable compromise for many applications involving nonlinear dynamics.

While the above specification remains within the $\mu$-linear framework $(RP)$, enabling the $\infty$-order variational analysis presented in Section~\ref{sec:dual}, for practical reasons, more sensitive $\mu$-nonlinear criteria could be utilized. One convenient option is the scalar characteristic:
\[
\text{Tr}(\mathbb{K} X) = \mathbb{E} \|X\|^2 - \big\|\mathbb{E} X\big\|^2 = \sum_{i=1}^n \mathbb{V}(X^i).
\]
If $X \sim \mu$, this value can be represented as:
\begin{equation}\label{Var}
\text{Tr}(\mathbb{K} X) = \frac{1}{2} \int_{\mathbb{R}^n} \mathrm{d}\mu(y) \int_{\mathbb{R}^n} \|x-y\|^2 \, \mathrm{d} \mu(x),
\end{equation}
which is a ``quadratic'' functional $\mathcal{P}(\mathbb{R}^n) \to \mathbb{R}$ of $\mu$. Fortunately, the latter can be reformulated as a linear form $\mathcal{P}(\mathbb{R}^n \times \mathbb{R}^n) \to \mathbb{R}$ by appropriately extending the state space. To achieve this, introduce a copy $Y$ of $X$ such that $Y_t \sim \mu_t$ for all $t \in I$, and $Y_0$ is independent of all previously addressed random variables, and consider an extended process $\bm{X} = (X, Y)$ driven by a $2m$-dimensional Brownian motion $\bm{W}$. Thus, $\bm{X}$ is defined on the product space $(\Omega\times \Omega, \mathcal{F} \otimes \mathcal{F},  \mathbf{P} \equiv \mathbb{P} \otimes \mathbb{P})$ as a solution of the SDE:
\begin{equation}
\begin{aligned}
\bm{X}_t &= \begin{bmatrix} X_0 \\ Y_0 \end{bmatrix} + \int_0^t \begin{bmatrix} f_s(X_s, u_s) \\ f_s(Y_s, u_s) \end{bmatrix} \, \mathrm{d} s + \int_0^t \begin{bmatrix} \sigma_s(X_s) & \mathbf{0}_{m\times m} \\ \mathbf{0}_{m\times m} & \sigma_s(Y_s) \end{bmatrix} \, \mathrm{d} \bm{W}_s, \label{SDE-2}
\end{aligned}
\end{equation}
assuming progressive measurability with respect to the filtration $t \mapsto \mathcal{F}^{\bm{W}, X_0, Y_0}_t$ generated by the process $\bm{W}$ and the random variables $X_0$, $Y_0$. Denoting the law of $\bm{X}$ by $\bm{\mu}$, the functional \eqref{Var} can now be expressed as:
\[
\text{Tr}(\mathbb{K} X) = \frac{1}{2} \int_{\mathbb{R}^n \times \mathbb{R}^n} \|x-y\|^2 \, \mathrm{d} \bm{\mu}(x, y).
\]
In this way, the corresponding optimal control problem is downshifted to the original statement $(P)$.

Finally, the same approach can be employed to ``linearize'' cost functionals of the ``$\mu$-polynomial'' structure involving terms such as:
\begin{equation}
\int_{\Omega} \ldots \int_{\Omega} \Psi(X_T(\omega^1), \ldots, X_T(\omega^{\rm q})) \, \mathrm{d} \mathbb{P}(\omega^1) \ldots \mathrm{d} \mathbb{P}(\omega^{\rm q}),\label{I-polyn}
\end{equation}
where $\Psi: \, (\mathbb{R}^n)^{\rm q} \to \mathbb{R}$. 

An obvious drawback of this transformation is the increase in dimensionality of the resulting optimal control problem $(P)$. A somewhat equivalent way to handle nonlinear costs of the form \eqref{I-polyn} is suggested in \cite{SPP-2023}, using the concept of higher-order adjoint to the law dynamics \eqref{PDE}. 

\section{The Approach}

In this section, we expose an approach to numerical solution of the problem $(RP)$, relying on an ``$\infty$-variation'' of $\mathcal J$ around a reference control.  

\subsection{Duality. Increment Formula}\label{sec:dual}

Fix $u \in \mathcal U$, and let $p=p[u] \colon I \times \R^n \to \R$ be a solution to the backward Kolmogorov equation
\begin{equation}
    \label{P-2}\{\partial_t  + \mathcal L_t\}p = 0, \quad p_T = \ell
\end{equation}
(the dependence on $u$ is dropped for brevity). It is well-known (see, e.g., \cite{Oksendal2010stochastic,Pham2021}) that, if this solution exists, it is unique and admits probabilistic representation via the  classical Feynman-Kac formula:
$
p_t(x) = \mathbb E\ell\left(X_{t, T}(x)\right),
$
in which $X_{s, t}(x)$ denotes a state of the SDE \eqref{SDE} at the moment $t \in [s, T]$ emanating from the deterministic position $X_s = x \in \R^n$ at the time moment $s \in [0, T)$. 

A direct computation with application of Ito's lemma leads to the pointwise equality
\[
\frac{d}{d t}\int p_t\d \mu_t = \int\big\{\partial_t + \mathcal L_t\big\}p_t\d \mu_t  = 0 \mbox{ a.a. }t\in I.
\]
The latter condition establishes the duality between $\mu$ and $p$.


Now, let $\bar u, u \in \mathcal{U}$ be arbitrary controls. The former stands for a reference (given) signal, while the latter is the target signal, intended to articulate the descent from $\bar u$. For simplicity, we use a bar to denote dependence on $\bar{u}$ and omit dependence on $u$, as in $\bar{p} \equiv p[\bar{u}]$ and $\mu \equiv \mu[u]$.

Consider the increment
\(
\Delta \mathcal I \doteq \mathcal I[u] - \mathcal I[\bar u]
\)
of the cost functional of the problem $(RP)$. By expressing it as 
\begin{align*}
\Delta {\mathcal J} & = \langle \mu_T - \bar \mu_T, \ell \rangle
 \doteq 
\langle \mu_T - \bar \mu_T, \bar p_T \rangle 
- 
\underbrace{\langle {\mu_0}-{\bar \mu_0}, \bar p_0 \rangle}_\text{$\equiv 0$}
\\& =\langle \mu_T, \bar p_T\rangle - \langle {\mu_0}, \bar p_0\rangle  =\int_I \frac{d}{d t}\langle \mu_t, \bar p_t\rangle\d t,
\end{align*}
computing 
\(
\frac{d}{dt} \langle {\mu_t}, \bar p_t\rangle = \langle {\mu_t}, \left\{\partial_t + \mathcal L\right\}\bar p\rangle,
\)
and evaluating
\[
\left\{\partial_t + \mathcal L\right\}\bar p 
=\left\{\mathcal L-\bar{\mathcal L}\right\}\bar p +  \underbrace{\left\{\partial_t + \bar{\mathcal{L}}\right\}\bar p}_\text{$\equiv 0$}
=
\left\{\mathcal L-\bar{\mathcal L}\right\}\bar{p} \doteq \nabla_x \bar p \cdot \left(f - \bar f\right),
\] 
we arrive at the desired representation
\begin{align}
\Delta \mathcal J & = \int_I \int_{\R^n}\left(\bar H_s(x, u_s) - \bar H_s(x, \bar u_s)\right) \d \mu_s(x) \d s.\label{increment-X-mu}
\end{align}
Here, 
\(
\bar H_s(x, \upsilon) \doteq H_s(x, \nabla_x \bar p_s(x), \upsilon)
\)
is a contraction of the usual Hamilton-Pontryagin function
$
H_s(x, \psi, \upsilon) \doteq \psi \cdot f_s(x, \upsilon)
$ 
to the gradient $\psi = \nabla_x  \bar p_s(x)$.





\subsection{Descent}\label{sec:decr}


Fix $\mu \in \mathcal{P}(\mathbb{R}^n)$ and $t \in I$, and take 
\begin{align}\label{fb}
\bar v_t[\mu] \in \arg\min_{\upsilon \in U}\int_{\R^n} \bar H_s(x, \upsilon) \d \mu(x).
\end{align}
Functions $\bar v\colon (t, \mu) \mapsto \bar v_t[\mu]$ serve as \emph{feedback controls} of the PDE \eqref{PDE}, whose substitution in the driving vector field $f$ results in the \emph{nonlocal} equation 
\begin{align}
      \partial_t\mu = \mathcal L_t^*(v_t[\mu]) \, \mu; \quad X_0 \sim \mu_0. \label{nlFPK}
\end{align}
Assume that \eqref{nlFPK} has a unique solution $\hat{\mu}=\hat\mu[\bar v]$ (in the above sense), and define 
\(
u(t) \doteq \bar v_t[\hat{\mu}_t].
\)
The use of this function as a target control in formula \eqref{increment-X-mu} evidently results in non-ascendancy in the cost functional: $\Delta \mathcal{I} \leq 0$. Repeated implementation of the described control-update rule $\bar u \to u$ gives rise to Algorithm~\ref{alg}. 
\begin{algorithm}
\caption{Descent method}
\label{alg}
\KwData{$\bar u \in \mathcal U$ (initial guess), $\varepsilon>0$ (tolerance)}
\KwResult{$\{u^k\}_{k \geq 0} \subset \mathcal U$ such that $\mathcal I[u^{k+1}] < \mathcal I[u^{k}]$}
$k \gets 0$;
$u^0 \gets \bar u$\;
\Repeat{$\mathcal I[u^{k-1}] - \mathcal I[u^{k}] < \varepsilon$}{
$p^{k} \gets p[u^k]$; $v^k_s[\mu]$ a solution of \eqref{fb} with $\bar u = u^k$\; $\mu^{k+1} \gets \hat \mu[v^k]$;
$u^{k+1} \gets v^k[\mu^{k+1}]$;
$k \gets k+1$\;
  }
\end{algorithm}

By construction, the sequence $\{\mathcal I^k \doteq \mathcal I[u^k]\}$ of cost values is monotone decreasing, ensuring convergence. However, it is important to note that $\bar v$ is generally discontinuous, rendering \eqref{nlFPK} ill-posed. 
This observation motivates giving preference to non-classical concepts of feedback solution, such as Krasovskii-Subbotin  (KS) ``constructive motions'', as discussed in \cite{SChPP-2022}. 

\subsection{Implementation}\label{sec:impl}

The conceptual algorithm outlined above relies on solving parabolic partial differential equations (PDEs) on $\mathbb{R}^n$, which can be challenging, especially for large values of $n$ (in fact, even when $n=4$). Implementing a discontinuous law-feedback control presents another difficulty. 

One potential solution is to revert to the original probabilistic framework and utilize Monte Carlo simulation of $\bar p$, along with an empirical approximation of $\mu$, combined with the KS sampling algorithm for sequential piecewise-constant synthesis of the descent control. 

Fixing $N, M, K \in \mathbb N$ and a partition $\pi = \{t_k\}$, $0=t_0 < t_1< \ldots < t_K=T$ of $I$, let us break down the probabilistic implementation of the Algorithm~\ref{alg}:

\vspace{0.1cm}

1) \emph{Approximation of the value $\bar p_t(x)$}. Using the Feynman-Kac formula,  we simulate the value $\bar p_{t}(x)$ at a point $(t, x) \in I \times \R^n$ using $N$ sample paths $\bar X^j$:
\(
    \bar p_{t}(x) = \mathbb E \psi(\bar X_{t, T}(x)) \approx \frac{1}{N} \sum_{j=1}^N\psi(\bar X^{j}_{t, T}(x)) \doteq  \bar p_{t}^N(x).
\)
This approximation is valid due to the uniform law of large numbers, ensuring convergence almost surely as $N \to \infty$ for any fixed $(t, x)$.

\vspace{0.1cm}

2) \emph{Approximation of the measure $\mu_t$}. Given a random state $X_t$, we approximate its law $\mu_t$ by empirical measures
\(
\mu^M_t \doteq \frac{1}{M} \sum_{j=1}^M \delta_{X^j_t},
\)
involving $M$ samples $X^j_t$ of $X_t$.\footnote{In general, $N \neq M$. In the above computations, we demonstrate a somewhat surprising vitality of this approach even when $M=1$.}

\vspace{0.1cm}

3)  \emph{Synthesis of a piecewise constant target control  ($k^{\text{th}}$ step)}. Assuming $u$ is calculated up to time $t_k \in \pi$, $k< K$, we update the constant value $w^{k}$ for the sub-interval $[t_{k}, t_{k+1}]$ using the KS sampling algorithm:

\vspace{0.1cm}

3.a) Draw $M$ sample paths $t \mapsto X^l_t = X_t[u]$, $t \in [0, t_k]$, and generate the population of points ${x_{lk} \equiv X^l_{t_k}}_{l=1}^M \subset \mathbb{R}^n$.  

\vspace{0.1cm}

3.b) Compute approximate values $\bar p_{t_k}^N(x_{lk})$, $l=1, \ldots, M$, and assign
\[
    u(t) \equiv w_k \in \arg\min_{\upsilon \in U}\frac{1}{M} \sum_{l=1}^M \bar H_{t_k}^{N}(x_{lk}, \upsilon).
\]

Note that the total computational cost of this iteration is proportional to $\#(N M K)$ solutions of the $n$-dimensional SDE \eqref{SDE}.

\subsection{Discussion}\label{sec:disc}

Monte Carlo approximation is known to be computationally expensive and prone to inaccuracies. However, when combined with the KS sampling algorithm, it offers a promising avenue for overcoming the curse of dimensionality associated with the numerical solution of PDEs \eqref{nlFPK} and \eqref{P-2}. This combination effectively reduces the cardinality of the computational mesh, with clear potential for parallelization. 

Importantly, the proposed ``approximate realization'' of Algorithm~\ref{alg}  \emph{loses the property of \textbf{monotone} descent}. Nevertheless, even with relatively coarse approximations of $\bar p$ and $\mu$, our approach exhibits surprising robustness in numerical experiments, demonstrating reasonably fast convergence ``on average''.



\subsection{Numerical Experiments}\label{experim}

    Consider a stochastic version of the Ermentrout–Kopell model of excitable neuron, known as ``Theta model'' \cite{Theta}.  The model represents a population of (identical and independent) excitable neurons characterized by their phase $X \in \mathbb S^1\doteq \R/2\pi\mathbb Z$ and the baseline current $Y$, assumed to be affected by the Brownian noise. Phases of all neurons are subject to common external excitations $u \in \R$. The dynamics of a sample neuron reads:
\begin{equation*}
    \dot X_t = (1-\cos X_t) +(1+\cos X_t)\left(Y_t+w(t, X_t, Y_t)\right), \quad
    d Y_t =  \sqrt{2\beta}\,d W_t.
\end{equation*}
The standard control problem consists in reaching a desired excitation pattern, which means that the neurons produce a spike at a predefined time instant $T$ with maximal probability. By accepting a convention that spiking is produced by passing through the phase $2\pi k$, $k \in \mathbb N$, our optimization problem is stated as:
\[
\ell(X_T) = (\sin(X_T))^{2\rm p} +  (\cos(X_T) - 1)^{2\rm p} \to \min, \quad \rm p \in \mathbb N.
\]

\begin{figure}\label{fig1}
  \centering
\includegraphics[width=1.\textwidth]{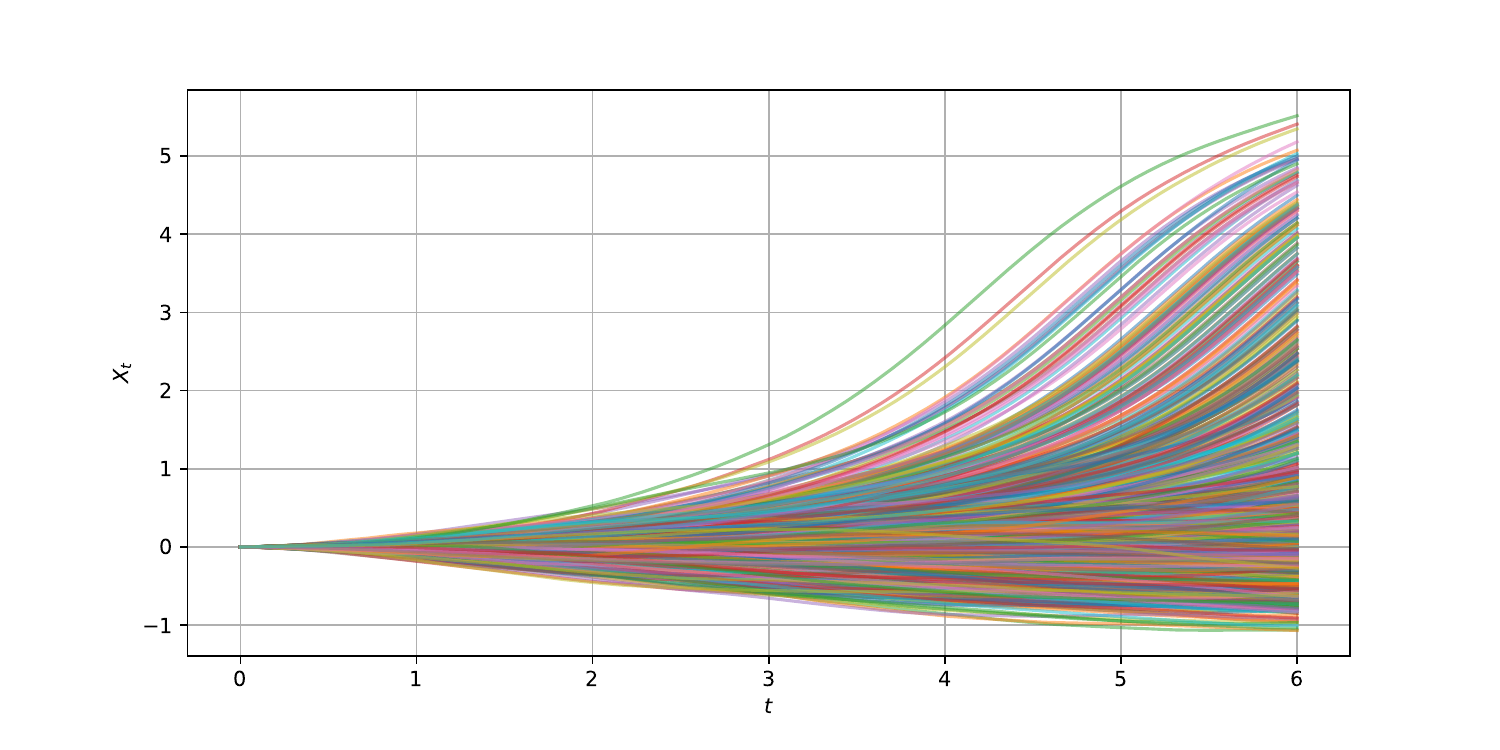}\vspace{-0.8cm}
  \includegraphics[width=1.\textwidth]{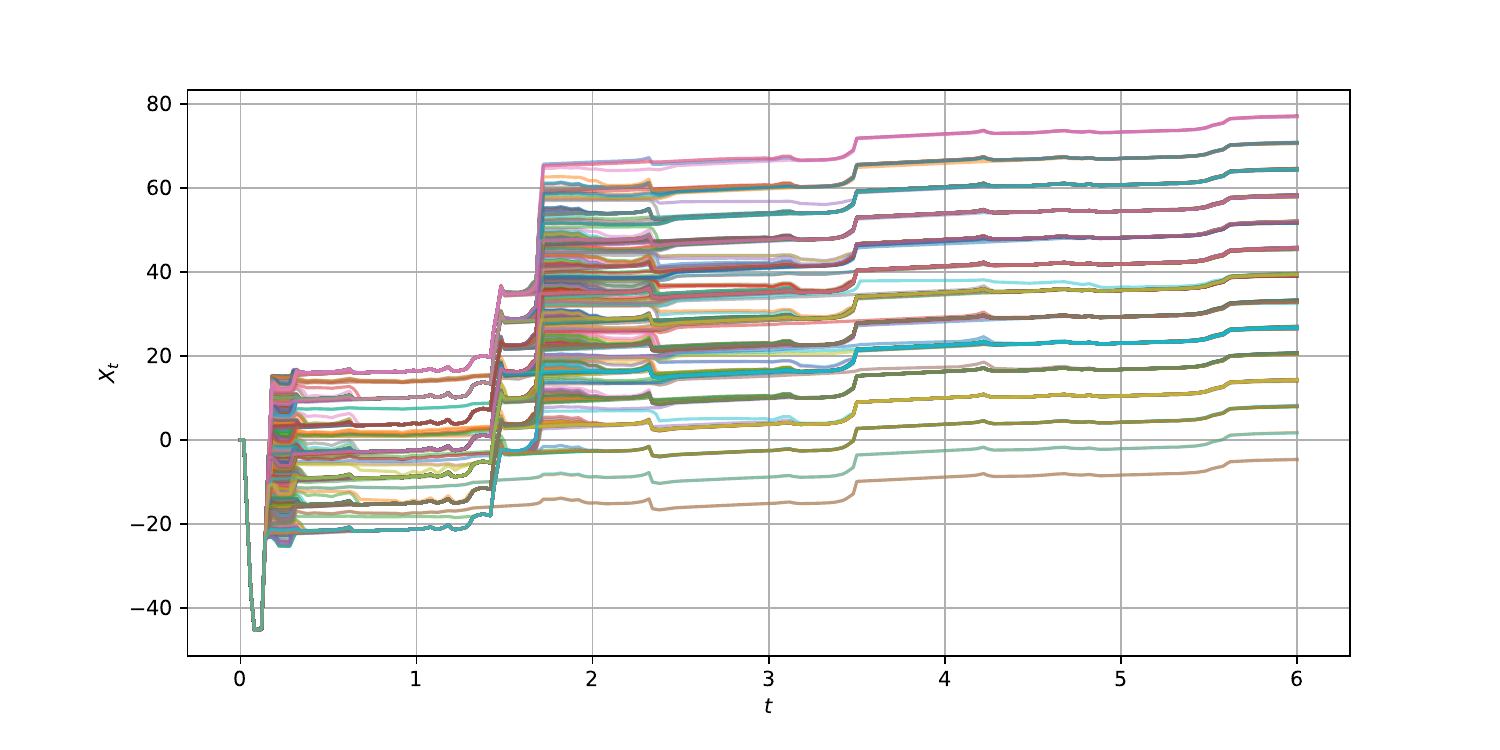}\vspace{-0.6cm}  
  \caption{Sample paths produced by the initial (upper) and learned (bottom) controls.}
\end{figure}

In our numerical experiments, we adhere to Markovian controls with the following predefined structure: \( w(t, x, y) = u_1(t) + u_2(t) y + u_3(t) \cos(x) + u_4(t) \sin(x) \), where \( u_j \) are measurable functions \( I \to \mathbb{R} \) subject to optimization. The computations are conducted for \( T = 6 \), \( \beta = 0.05 \), \( \text{p} = 1, 2\), and \( U = \mathbb{R} \), assuming quadratic penalties \( \sum_j u_j^2 \) for minimization in \eqref{fb}. The control \( u = (u_1, \ldots, u_4) \) is learned using Algorithm 1 with the approximation discussed in \S~2.1 with parameters \( N = 100 \), \( M = 1 \), and \( K = 20 \) per unit time, starting with \( u^0 = (0, 0, 0, 0) \). The learned control is tested using \( N = 1000 \) samples, and the averaged performance \( \check{\mathcal{I}} \) is computed. For \(\text{p} = 1\), optimization is achieved in 3 iterations with $\check{\mathcal{I}}^0 \approx 2.39$ and $\check{\mathcal{I}}^3 \approx 0.02$ 
Fig.~\ref{fig1} exposes the uncontrolled and controlled populations \( t \mapsto X_t \) for the case \(\text{p} = 1\). Interestingly, the effect of ``quantization'' (steering different clusters of the population to different equivalent phases \( 2\pi k \), \( k \in \mathbb{N} \)) is observed. As expected, an even stronger denoising effect is achieved for the higher-order statistical tracking \(\text{p} = 2\).

\begin{credits}
\subsubsection{\ackname} 

The authors acknowledge the financial support of the Foundation for Science and
Technology (FCT, Portugal) in the framework of the Associated Laboratory ARISE (LA/P/0112/2020), R\&D Unit SYSTEC (base UIDB/00147/2020 and programmatic UIDP/00147/2020 funds), and project RELIABLE (PTDC/EEI-AUT/3522/2020). A part of the computations was carried out on the OBLIVION Supercomputer (Évora University) under FCT computational project 2022.15706.CPCA. RC and MS acknowledge personal financial support by the FCT with DOI refs.: 10.54499/CEECINST/00010/2021/CP1770/CT0006 and 
10.54499/CEECINST/00049/ 2018/CP1524/CT0006, respectively. 

\end{credits}

\bibliographystyle{splncs}
\bibliography{Staritsyn-full}
%




\end{document}